\documentclass[11pt]{amsart}

\usepackage{amssymb}
\usepackage{enumitem}
\usepackage{hyperref}
\usepackage{tikz}

\usepackage{fullpage}

\newtheorem{theorem}{Theorem}[section]
\newtheorem{lemma}[theorem]{Lemma}

\newtheorem{alphtheorem}{Theorem}
\newtheorem{alphcorollary}[alphtheorem]{Corollary}

\theoremstyle{definition}

\newtheorem{remark}[theorem]{Remark}

\newtheorem{definition}[theorem]{Definition}

\numberwithin{equation}{section}
\numberwithin{figure}{section}

\DeclareMathOperator{\Pic}{Pic}

\DeclareMathOperator{\Hom}{Hom}

\DeclareMathOperator{\ev}{ev}
\DeclareMathOperator{\supp}{supp}

\DeclareMathOperator{\Hilb}{Hilb}

\DeclareMathOperator{\nc}{nc}
\let\epsilon\varepsilon

\newcommand*\ratmap{\mathbin{\tikz [baseline=0ex,-latex, dashed, ->] \draw [densely dashed] (0em,0.58ex) -- (1.3em,0.58ex);}}

\title[Maps into abelian varieties]{On the finiteness of maps into simple abelian varieties satisfying certain tangency conditions}

\author{Finn Bartsch}
\address{Finn Bartsch \\
IMAPP Radboud University Nijmegen \\
PO Box 9010, 6500GL \\
Nijmegen, The Netherlands\\}
\email{f.bartsch@math.ru.nl}


\begin{document}

\begin{abstract}
We show that given a simple abelian variety $A$ and a normal variety $V$ defined over a finitely generated field $K$ of characteristic zero, the set of non-constant morphisms $V \to A$ satisfying certain tangency conditions imposed by a Campana orbifold divisor $\Delta$ on $A$ is finite.
To do so, we study the geometry of the scheme $\underline{\Hom}^{\nc}(C, (A, \Delta))$ parametrizing such morphisms from a smooth curve $C$ and show that it admits a quasi-finite non-dominant morphism to $A$.
\end{abstract}

\maketitle
\thispagestyle{empty}

\section{Introduction}

The purpose of this note is to prove a peculiar finiteness result for maps into simple abelian varieties satisfying certain tangency conditions.
To state the result, we employ the notion of a \emph{Campana orbifold} (also appearing in the literature under the name \emph{C-pair}).
These objects are pairs $(X, \Delta)$ consisting of a smooth variety $X$ (we will always assume varieties to be geometrically integral) and a $\mathbb{Q}$-divisor on $X$ whose coefficients are of the form $1-\frac{1}{m}$ with $m \in \mathbb{N}_{\geq 1} \cup \{ \infty \}$, where we interpret $\frac{1}{\infty} = 0$.
If $\Delta = \sum_i (1-\frac{1}{m_i}) D_i$ is a Campana orbifold $\mathbb{Q}$-divisor, we refer to $m_i$ as the \emph{multiplicity} of $D_i$ in $\Delta$.
A Campana orbifold $(X, \Delta)$ is said to be \emph{smooth} if $\supp \Delta$ is a simple normal crossings divisor.
If $V$ is a normal variety and $(X, \Delta)$ is a Campana orbifold, then a morphism $f \colon V \to (X, \Delta)$ is a morphism of varieties $f \colon V \to X$ satisfying $f(V) \subsetneq \supp \Delta$ such that for every prime divisor $D \subseteq \supp \Delta$ and every prime divisor $E \subseteq V$ appearing in $f^*D$, the coefficient of $E$ in $f^*D$ is at least the multiplicity of $D$ in $\Delta$.
Note that if $\Delta = 0$ is the empty divisor, then a morphism $V \to (X,0)$ is just a morphism of varieties $V \to X$.
This allows us to view every smooth variety as a Campana orbifold. 
A \emph{near-map} $V \ratmap (X, \Delta)$ is a rational map $V \ratmap X$ for which there is a closed subset $Z \subsetneq V$ of codimension at least two such that $V \setminus Z \ratmap (X, \Delta)$ is a morphism.

\begin{alphtheorem}\label{main_thm}
Let $K$ be a field which is finitely generated over $\mathbb{Q}$.
Let $A$ be an abelian variety over $K$ which is simple over $\overline{K}$.
Let $\Delta$ be a nonzero Campana orbifold $\mathbb{Q}$-divisor on $A$ whose support is a simple normal crossings divisor.
Then, for every normal variety $V$ over $K$, the set of nonconstant near-maps $V \ratmap (A, \Delta)$ defined over $K$ is finite.
\end{alphtheorem}

To prevent any possible confusion, let us stress that ``Campana orbifolds'' are not really orbifolds in the classical sense of the word (i.e.\ they are not finite type separated smooth DM-stacks with trivial generic stabilizer).
Nonetheless, Theorem \ref{main_thm} readily implies the following statement about maps into root stacks.
Essentially, this is because a morphism $V \to \sqrt[m]{D/A}$ into the root stack $\sqrt[m]{D/A}$ induces a Campana orbifold morphism $V \to (A, \Delta)$ with $\Delta := (1-\frac{1}{m})D$ unless it set-theoretically factors through $D$ (see \cite[Section~10.3]{OlssonBook} for the definition of root stacks).
However, note that the converse is not true, i.e.\ most Campana orbifold morphisms $V \to (A, \Delta)$ will not induce a morphism $V \to \sqrt[m]{D/A}$.

\begin{alphcorollary}\label{cor_rootstacks}
Let $K$ be a field which is finitely generated over $\mathbb{Q}$.
Let $A$ be an abelian variety over $K$ which is simple over $\overline{K}$.
Let $D \subseteq A$ be a smooth irreducible divisor and let $m \geq 2$ be an integer.
Then, for every variety $V$ over $K$, the set of nonconstant morphisms into the root stack $V \to \sqrt[m]{D/A}$ defined over $K$ is finite.
\end{alphcorollary}

In our proof of Theorem \ref{main_thm}, we first reduce to the case that $V = C$ is one-dimensional by repeatedly intersecting $V$ with hyperplanes.
The analysis of the one-dimensional case then crucially relies on the fact that the set of morphisms from the smooth curve $C$ into $(A, \Delta)$ naturally comes equipped with a scheme structure.
We proceed by examining the geometry of this morphism scheme.
More precisely, we will make use of the following recently established result.
We say that a Campana orbifold $(X, \Delta)$ is \emph{projective} if $X$ is so.

\begin{theorem}\label{orbifold_hom_scheme}
Let $k$ be a field of characteristic zero.
Let $C$ be a smooth quasi-projective curve over $k$ with smooth projective model $\overline{C}$.
Let $(X, \Delta)$ be a projective Campana orbifold.
Then there is a closed subscheme
\[ \underline{\Hom}_k(C, (X, \Delta)) \subseteq \underline{\Hom}_k(\overline{C}, X) \setminus \underline{\Hom}_k(\overline{C}, \supp \Delta) \]
satisfying the following properties:
\begin{enumerate}
\item $\underline{\Hom}(C, (X, \Delta))(k)$ is the set of orbifold morphisms $C \to (X, \Delta)$ defined over $k$.
\item If $k \subseteq l$ is an algebraic extension, then $\underline{\Hom}_l(C_l, (X_l, \Delta_l)) = (\underline{\Hom}_k(C, (X, \Delta)))_l$.
\item For every irreducible component $H \subseteq \underline{\Hom}_k(C, (X, \Delta))$ with normalization $H'$, the evaluation morphism $\ev \colon C \times H' \to (X, \Delta)$ is an orbifold morphism.
\end{enumerate}
\end{theorem}
\begin{proof}
See \cite[Theorem~3.9~and~Corollary~3.10]{BJRWeaklySpecial}.
\end{proof}

Studying the geometry of this $\Hom$-scheme, we obtain the following result, from which we then deduce Theorem \ref{main_thm} via an application of Faltings's finiteness theorem for rational points on subvarieties of abelian varieties.
We stress that this next statement is a purely geometric statement, i.e.\ it is a statement over $\overline{K}$ -- the arithmetic input comes entirely from Faltings's theorem.

\begin{alphtheorem}\label{main_thm2}
Let $k$ be an algebraically closed field of characteristic zero.
Let $A$ be a simple abelian variety over $k$ and let $\Delta$ be a nonzero Campana orbifold divisor on $A$ whose support is a simple normal crossings divisor.
Let $C$ be a smooth quasi-projective curve over $k$.
Then, the scheme $\underline{\Hom}_k^{\nc}(C,(A, \Delta))$ parametrizing non-constant morphisms $C \to (A, \Delta)$ has only finitely many irreducible components, each of which admits a non-dominant map to $A$ with finite fibers.
\end{alphtheorem}

Let us briefly comment on the necessity of the assumptions in Theorems \ref{main_thm} and \ref{main_thm2}.

\begin{remark}
Let us give an example showing that the assumption that $A$ is simple in Theorem \ref{main_thm} is necessary.
For this, let $A := A' \times A''$ where $A'$ and $A''$ are positive-dimensional abelian varieties over $K$ such that $A''(K)$ is infinite.
Let $D'$, respectively $D''$, denote any nonzero Campana orbifold divisor on $A'$, respectively $A''$, and consider $D := D' \times A'' + A' \times D''$ on $A' \times A''$.
Let $V$ be any normal variety admitting a nonconstant morphism $V \to (A', D')$.
Then, the infinitely many $K$-rational points of $A''$ lead to infinitely many horizontal embeddings $A' \to A' \times A'' = A$ defined over $K$.
Composing these with the nonconstant morphism $V \to (A', D')$ gives rise to infinitely many orbifold morphisms $V \to (A, D)$.
Similar examples can be given if $A$ is merely isogeneous to a product.
Observe also that this issue cannot be easily circumvented by strengthening the assumptions on the Campana orbifold divisor, as $D$ might very well be ample.
\end{remark}

\begin{remark}
Let us give an example showing that the components of the $\Hom$-scheme appearing in Theorem \ref{main_thm2} can be positive-dimensional.
This also shows that the assumption that $K$ is finitely generated over $\mathbb{Q}$ is necessary in Theorem \ref{main_thm}.
Let $C$ be a curve of genus $2$ whose Jacobian is simple and let $\sigma \colon C \to C$ denote its hyperelliptic involution.
Let $\Theta$ be the image of the natural map $C \to \Pic^1(C)$.
Then, for all $c, d \in C$, we can consider the translate $\Theta + [c] - [d] \subseteq \Pic^1(C)$.
It is distinct from $\Theta$ as long as $c \neq d$.
One easily checks that the curve $\Theta$ and its translate $\Theta + [c] - [d]$ intersect exactly in the two points $[d]$ and $[\sigma(c)]$.
Thus, for every $c \in C$ such that $c \neq \sigma(c)$, we have that $\Theta$ is distinct from its translate $\Theta + [c] - [\sigma(c)]$, and intersects it exactly in the point $[\sigma(c)]$ with intersection multiplicity $2$.
This gives rise to a positive-dimensional family of morphisms $C \to (\Pic^1(C), \frac{1}{2} \Theta)$, indexed by the open subset of $C$ consisting of the points which are not fixed points of $\sigma$.
\end{remark}

\subsection{Relation to other work}

Our Theorem~\ref{main_thm} can be viewed as saying that, for $A$ a simple abelian variety and $\Delta$ a nontrivial Campana orbifold divisor with simple normal crossings support, the Campana orbifold $(A, \Delta)$ is \emph{$p$-Mordellic} in the sense of \cite[Definition~1.2]{EJRIntermediate} for every integer $p \geq 1$, where we just naively adapt the definition given there to the orbifold context.
In particular, since the smooth projective orbifold pair $(A, \Delta)$ is of general type, this can be seen as a partial confirmation of Conjecture~1.4 in loc.\ cit.\ .

In \cite[Theorem~1.6]{JavMRL}, Javanpeykar proves that if $K$ is finitely generated over $\mathbb{Q}$ and $A$ is an abelian variety which is simple over $\overline{K}$, then for every normal variety $X$ admitting a finite surjective ramified morphism $X \to A$ and every normal variety $V$, there are only finitely many non-constant rational maps $V \ratmap X$.
If we additionally assume that the morphism $X \to A$ is Galois and its branch locus has a smooth component, then we can rederive this result from our Theorem~\ref{main_thm}.
Indeed, if $D \subseteq A$ is a smooth component of the branch locus of $X \to A$ and $m$ is the order of ramification above it, then every rational map $V \ratmap X$ induces a Campana orbifold near-map $V \ratmap (A, (1-\frac{1}{m})D)$.
Hence the finiteness of the maps into $(A, (1-\frac{1}{m}) D)$ implies the finiteness of the maps into $X$.
We note that while the smoothness assumption on $D$ is probably only a technical artifact of our proof, the fact that this line of reasoning only seems to work for Galois covers $X \to A$ appears to be a more structural problem.

\subsection{Acknowledgements}

I thank Ariyan Javanpeykar for constant support and many helpful discussions.

\section{Proof of the main results}

Let us now come to the proof of Theorem \ref{main_thm2}.
It will be convenient to formulate our proof using various notions of hyperbolicity which we now introduce.

\begin{definition}\label{def:generaltype}
Let $k$ be a field.
Let $(X, \Delta)$ be a smooth projective Campana orbifold over $k$.
We say that $(X, \Delta)$ is of \emph{general type} if the canonical divisor of the pair $(X, \Delta)$, i.e.\ $K_X + \Delta$, is big.
\end{definition}

The following notion was studied by Javanpeykar--Kamenova under the name ``$(1,1)$-boun\-ded\-ness'' \cite[Definition~4.2]{JavanpeykarKamenova}.
Their definition naturally extends to the setting of Campana orbifolds.

\begin{definition}\label{def:geomhyp}
Let $(X, \Delta)$ be a projective Campana orbifold over an algebraically closed field $k$ of characteristic zero and let $Z \subseteq X$ be a closed subset.
We say that $(X, \Delta)$ is \emph{geometrically hyperbolic modulo $Z$} if for every closed point $x \in X \setminus (Z \cup \supp \Delta)$, every smooth quasi-projective curve $C$ over $k$ and every closed point $c \in C$, the set of Campana orbifold morphisms $C \to (X, \Delta)$ mapping $c$ to $x$ is finite.
We say that $(X, \Delta)$ is \emph{geometrically hyperbolic} if $(X, \Delta)$ is geometrically hyperbolic modulo the empty set.
\end{definition}

A Campana orbifold $(X, \Delta)$ being geometrically hyperbolic forces the $\Hom$-schemes of morphisms from a curve into $(X, \Delta)$ to be low-dimensional, even for ``unpointed'' morphisms.
We use the superscript $~^{\nc}$ to denote the part of the $\Hom$-scheme which parametrizes nonconstant morphisms. 

\begin{lemma}\label{geom_hyp_dimension_estimate}
Let $k$ be an algebraically closed field of characteristic zero.
Let $(X, \Delta)$ be a smooth projective Campana orbifold of general type over $k$ and let $Z \subseteq X$ be a closed subset such that $(X, \Delta)$ is geometrically hyperbolic modulo $Z$.
Then, for every smooth quasi-projective curve $C$ over $k$, we have
\[ \dim \underline{\Hom}_k^{\nc}(C, (X, \Delta)) \setminus \underline{\Hom}_k^{\nc}(C, Z) < \dim X. \]
\end{lemma}
\begin{proof}
Let $H \subseteq \underline{\Hom}_k^{\nc}(C, (X, \Delta))$ be any irreducible component not contained in $\underline{\Hom}_k^{\nc}(C,Z)$. 
Assume for a contraction that $\dim H \geq \dim X$.
Observe that for any fixed point $c \in C$, the evaluation morphism $\ev_c \colon H \to X$ has finite fibers over $X \setminus (Z \cup \supp \Delta)$:
If there was a point $x \in X \setminus (Z \cup \supp \Delta)$ such that the fiber of $\ev_c$ over $x$ was infinite, we would have found infinitely many morphisms $C \to (X, \Delta)$ mapping $c$ to $x$, contradicting the geometric hyperbolicity of $(X, \Delta)$ modulo $Z$.
As the morphism $\ev_c$ does not factor over $Z \cup \supp \Delta$ for almost every $c \in C$ by our choice of $H$, we see that for almost every $c \in C$, the morphism $\ev_c \colon H \to X$ is dominant.
Moreover, by Theorem \ref{orbifold_hom_scheme}, we know that the morphisms $\ev_c \colon H \to (X, \Delta)$ are orbifold.
Hence, we have found infinitely many dominant orbifold morphisms to the Campana orbifold $(X, \Delta)$.
However, as $(X, \Delta)$ is of general type, by \cite[Theorem~1.1]{BJKobayashiOchiai} there are only finitely many dominant orbifold morphisms $H \to (X, \Delta)$.
This means that infinitely many of the morphisms $\ev_c$ must coincide, which by continuity means that all $\ev_c$ must be equal to each other.
But this contradicts the assumption that $H$ contains only nonconstant morphisms, concluding the proof.
\end{proof}

Another notion of hyperbolicity we will use is Demailly's notion of algebraic hyperbolicity, suitably adapted to the setting of Campana orbifolds.

\begin{definition}\label{def:alghyp}
Let $(X, \Delta)$ be a projective Campana orbifold over an algebraically closed field $k$ of characteristic zero and let $Z \subseteq X$ be a closed subset.
We say that $(X, \Delta)$ is \emph{algebraically hyperbolic modulo $Z$} if for one (then every) ample line bundle $\mathcal{L}$ on $X$, there is a constant $\epsilon > 0$ such that for every smooth quasi-projective curve $C$ with smooth projective model $\overline{C}$ of geometric genus $g_{\overline{C}}$ and every non-constant morphism $f \colon C \to (X, \Delta)$ with extension $\overline{f} \colon \overline{C} \to X$ satisfying $f(C) \nsubseteq Z$, we have that
\[ \deg_{\overline{C}} \overline{f}^*\mathcal{L} \leq \epsilon (2g_{\overline{C}} - 2 + \#(\overline{C} \setminus C)). \]
We say that $(X, \Delta)$ is \emph{algebraically hyperbolic} if $(X, \Delta)$ is algebraically hyperbolic modulo the empty set.
\end{definition}

The following Lemma is well-known, for the convenience of the reader we give a proof.

\begin{lemma}\label{alg_hyp_bounded}
Let $k$ be an algebraically closed field of characteristic zero.
Let $(X, \Delta)$ be a projective Campana orbifold over $k$ and let $Z \subseteq X$ be a closed subset such that $(X, \Delta)$ is algebraically hyperbolic modulo $Z$.
Then, for every smooth quasi-projective curve $C$ over $k$, the scheme $\underline{\Hom}_k(C, (X, \Delta)) \setminus \underline{\Hom}_k(C, Z)$ has finitely many irreducible components.
\end{lemma}
\begin{proof}
Let $\overline{C}$ be the smooth projective model of $C$.
The scheme $\underline{\Hom}_k(C, (X, \Delta))$ is by construction a subscheme of the Hilbert scheme $\Hilb_k(\overline{C} \times X)$.
Hence, it suffices to show that the graphs of the morphisms $f \colon \overline{C} \to X$ parametrized by $\underline{\Hom}_k(C, (X, \Delta)) \setminus \underline{\Hom}_k(C, Z)$, viewed as closed subschemes of $\overline{C} \times X$, have only finitely many Hilbert polynomials with respect to any fixed ample line bundle on $\overline{C} \times X$.  
Let $\mathcal{L}$ be an ample line bundle on $X$ and let $c \in \overline{C}$ be any closed point.
Then, the line bundle $\widetilde{\mathcal{L}} := \mathcal{O}_{\overline{C}}(c) \boxtimes \mathcal{L}$ on $\overline{C} \times X$ is ample.
If $f \colon \overline{C} \to X$ is any morphism and $\Gamma_f \subseteq \overline{C} \times X$ is its graph, then the restriction of $\widetilde{\mathcal{L}}$ to $C \cong \Gamma_f$ is given by $\mathcal{O}_{\overline{C}}(c) \otimes f^*\mathcal{L}$.
Consequently, the Hilbert polynomial of $\Gamma_f$ with respect to $\widetilde{\mathcal{L}}$ is $n \mapsto (1-g) + (1+\deg_{\overline{C}}(f^*\mathcal{L}))n$, where $g$ is the genus of $\overline{C}$.
By assumption, $\deg_{\overline{C}}(f^*\mathcal{L})$, which is always a positive integer, is bounded from above as $f$ ranges over those morphisms $\overline{C} \to X$ which represent orbifold morphisms $C \to (X, \Delta)$ whose image is not contained in $Z$.
This finishes the proof.
\end{proof}

The next Lemma is essentially a consequence of Yamanoi's work on the Green--Griffiths--Lang conjecture for varieties with maximal Albanese dimension \cite{YamanoiMaxAlbanese}.

\begin{lemma}\label{a_delta_hyperbolic}
Let $k$ be an algebraically closed field of characteristic zero.
Let $A$ be a simple abelian variety over $k$ and let $\Delta$ be a nonzero Campana orbifold divisor on $A$ whose support has simple normal crossings.
Then $(A, \Delta)$ is of general type, geometrically hyperbolic, and algebraically hyperbolic.
\end{lemma}
\begin{proof}
As $A$ is simple, every effective divisor on $A$ is ample (see e.g.\ \cite[Section~6, Application~1]{MumfordAV}).
This shows that $(A, \Delta)$ is of general type.

By \cite[Theorems~3.21~and~3.22]{JRGeomSpecial} (which heavily rely on \cite{YamanoiMaxAlbanese}), there is a proper closed subset $E \subsetneq A$ such that $(A, \Delta)$ is both algebraically and geometrically hyperbolic modulo $E$.
Since $A$ is simple, every proper closed subvariety $Z \subsetneq A$ is algebraically hyperbolic (either combine \cite[Theorem~2]{KawamataBloch} and \cite[Theorem~2.1]{DemaillyAlgHyp} or use the more general \cite[Corollary~1.(3)]{YamanoiMaxAlbanese}) and hence geometrically hyperbolic \cite[Theorem~10.8~and~Corollary~11.5]{JavanpeykarSurvey}.
In particular, every irreducible component of $E$ is both algebraically and geometrically hyperbolic.
This implies the claim.
\end{proof}

We can now give the proof of the results announced in the introduction.

\begin{proof}[Proof of Theorem \ref{main_thm2}]
By Lemma \ref{a_delta_hyperbolic}, the Campana orbifold $(A, \Delta)$ is of general type, geometrically hyperbolic, and algebraically hyperbolic.
Consequently, by Lemma \ref{alg_hyp_bounded}, we see that $\underline{\Hom}_k^{\nc}(C, (A,\Delta))$ has only finitely many irreducible components.
Let $H$ be one of these irreducible components, endowed with the reduced scheme structure.
By Lemma \ref{geom_hyp_dimension_estimate}, we have $\dim H < \dim A$.
Now, for a closed point $c \in C$, we can consider the restriction of the evaluation map $\ev \colon C \times H \to A$ to $\{ c \} \times H$.
Let $\ev_c \colon H \to A$ be the map thus obtained.
As $(A, \Delta)$ is geometrically hyperbolic, the map $\ev_c$ has finite fibers and it cannot be dominant for dimension reasons.
\end{proof}

\begin{remark}
While the statement of Theorem \ref{main_thm2} is enough to deduce Theorem \ref{main_thm}, we can prove slightly more.
Namely, if we again let $H$ be an irreducible component of $\underline{\Hom}_k^{\nc}(C, (A,\Delta))$ and let $\overline{H}$ denote the closure of $H$ in $\underline{\Hom}_k(C,A)$, the scheme $\overline{H}$ is a projective variety admiting a finite morphism to $A$.
To see this, first note that since $A$ contains no rational curves, the irreducible components of $\underline{\Hom}_k(C,A)$ are projective by \cite[Lemma~3.5~and~Proposition~3.9]{JavanpeykarKamenova}.
Thus, $\overline{H}$ is a projective variety.
Now let $c \in C$ be a closed point for which the evaluation map $\ev_c \colon \overline{H} \to A$ does not factor over $\supp \Delta$.
Then, by the previous, if $F \subseteq \overline{H}$ is a positive-dimensional subvariety contracted by $\ev_c$, we must have $F \subseteq \overline{H} \setminus H$.
By Theorem \ref{orbifold_hom_scheme}, the points of $\overline{H} \setminus H$ parametrize morphisms $C \to \supp \Delta$.
As $A$ is a simple abelian variety, the irreducible components of $\supp \Delta$ are geometrically hyperbolic by the aforementioned \cite[Corollary 1.(3)]{YamanoiMaxAlbanese} and \cite[Theorem~10.8~and~Corollary~11.5]{JavanpeykarSurvey}.
This geometric hyperbolicity contradicts the positive-dimensionality of $F$.
Hence $\ev_c$ is a proper morphism with finite fibers and thus is finite, as claimed.
\end{remark}

\begin{proof}[Proof of Theorem \ref{main_thm}]
We first reduce to the case that $V = C$ is a smooth quasi-projective curve.
For this, assume that $V$ was a variety admitting infinitely many pairwise distinct near-maps $(f_i \colon V \ratmap (A, \Delta))_{i \in \mathbb{N}}$ defined over $K$.
Replacing $V$ by an open subset if necessary, we may assume that $V$ is smooth and quasi-projective.
Fix an immersion $V \subseteq \mathbb{P}^r$ into some projective space.
After choosing an embedding $K \subseteq \mathbb{C}$, the infinitely many base changed near-maps $(f_{i,\mathbb{C}} \colon V_\mathbb{C} \ratmap (A_\mathbb{C}, \Delta_\mathbb{C}))$ are still pairwise distinct.
Now let $W \subseteq V_\mathbb{C}$ be a very general hyperplane section.
As long as $V$ is at least two-dimensional, $W$ will be (geometrically) integral, hence a variety.
Since $W$ is very general, we have that $f_i(W) \nsubseteq \Delta_\mathbb{C}$ for every $i$.
Thus, the restrictions $f_i|_W$ define Campana orbifold near-maps $W \ratmap (A_\mathbb{C}, \Delta_\mathbb{C})$ (see \cite[Corollary~2.7]{BJKobayashiOchiai}).
Moreover, $W$ being very general implies that the restrictions $f_i|_W$ are still pairwise distinct.
While $W$ may not be definable over $K$, there is a finitely generated extension $K \subseteq L$ such that $W$ can be defined over $L$.
Thus, at the cost of replacing $K$ by a finitely generated extension, we have reduced the dimension of $V$ by one.
Iterating this process, we hence may assume that $V = C$ is one-dimensional, as desired.
In particular, the $f_i$ now all define genuine morphisms $C \to (A, \Delta)$.

Now consider the orbifold $\Hom$-scheme $\underline{\Hom}_K^{\nc}(C, (A,\Delta))$.
By applying Theorem \ref{main_thm2} to its base change $(\underline{\Hom}_K^{\nc}(C, (A,\Delta)))_{\overline{K}}$, we see that it has only finitely many irreducible components.
Moreover, we see that, possibly after replacing $K$ by a finite extension, each of these irreducible components admits a non-dominant morphism to $A$ with finite fibers.
Thus, by applying Faltings's theorem on $K$-rational points of closed subvarieties of abelian varieties \cite[Theorem~5.1]{FaltingsGeneralSLang} to the image of these irreducible components in $A$, we see that each of them has only finitely many points over $K$.
This concludes the proof.
\end{proof}

\begin{proof}[Proof of Corollary \ref{cor_rootstacks}]
Replacing $V$ by its normalization if necessary, we may assume that $V$ is normal.
If $V \to \sqrt[m]{D/A}$ is a nonconstant morphism, we can look at the induced morphism $f \colon V \to A$ to the coarse space.
Unless $f$ factors over $D$, the pullback $f^* D$ is well-defined as a Weil divisor on $V$.
Because $f$ was induced from a morphism to the root stack $\sqrt[m]{D/A}$, we see that the coefficients of the divisor $f^* D$ are all divisible by $m$.
In particular, $f$ induces a Campana orbifold morphism $V \to (A, \Delta)$ with $\Delta := (1-\frac{1}{m}) D$.
Since we know that the set of morphisms $V \to (A, \Delta)$ is finite by Theorem~\ref{main_thm}, we only have to prove finiteness of those morphisms which factor over $D$.
But this finiteness follows directly from the finiteness of the $K(V)$-rational points of $D$, which is again a consequence of Faltings's theorem \cite[Theorem~5.1]{FaltingsGeneralSLang}.
\end{proof}

\bibliographystyle{alpha}
\bibliography{defranchis}{}
\end{document}